   \documentclass[11pt]{amsart}
 \usepackage{fancyhdr}
  \usepackage{times}
\usepackage[T1]{fontenc}
\usepackage{latexsym}
\usepackage[dvips]{graphics}
\usepackage{epsfig}
 \usepackage{hyperref, amsmath, amsthm, amsfonts, flafter,epsf}
\usepackage{amssymb,amscd}
\input amssym.def
\input amssym.tex
\usepackage{color}
        \usepackage[parfill]{parskip}    
 \usepackage{mathrsfs} 
\usepackage{epstopdf}

\newtheorem{thm}{Theorem}

\newtheorem{thm*}{Theorem}
\newtheorem*{hyp 1}{Weak Mertens Hypothesis}
\newtheorem*{hyp 2}{Linear Independence Hypothesis}
\newtheorem*{con*}{Conjecture}
\newtheorem*{lem*}{Lemma}

\theoremstyle{definition}

\newcommand{\ds}{\displaystyle}

\renewcommand{\d}{\delta}
\renewcommand{\b}{\beta}
\newcommand{\g}{\gamma}
\newcommand{\e}{\epsilon}

\renewcommand{\r}{\rho}
\newcommand{\Z}{\zeta}
\newcommand{\s}{\sigma}
\newcommand{\hh}{\tfrac12}

\newcommand\be{\begin{equation}}
\newcommand\ee{\end{equation}}
\newcommand\bea{\begin{eqnarray}}
\newcommand\eea{\end{eqnarray}}
\newcommand\bi{\begin{itemize}}
\newcommand\ei{\end{itemize}}
\newcommand\ben{\begin{enumerate}}
\newcommand\een{\end{enumerate}}

\begin{document}

 \title[ ]{Remarks on a Formula of Ramanujan}
\author{Andr\'es Chirre}
\author{Steven M. Gonek}
 
 
 \address{Departamento de Ciencias - Sección Matemáticas, Pontificia Universidad Católica del Perú, Lima, Perú}
\email{cchirre@pucp.edu.pe}

\address{Department of Mathematics, University of Rochester, Rochester, NY 14627, USA}
\email{gonek@math.rochester.edu \quad  }



\maketitle

\begin{abstract}
Assuming an averaged form of Mertens' conjecture and that the ordinates of the non-trivial zeros of the Riemann zeta function are linearly independent over the rationals, we analyze the finer structure of the terms in a well-known formula of Ramanujan.
\end{abstract}



\section{The formula of Ramanujan}\label{Ram form}
Let  $\mu(n)$ be the M\"obius function and set
\be\notag
F(b) = \sum_{n=1}^{\infty} \frac{\mu(n)}{n} e^{-(b/n)^2} .
\ee
In their paper ``Contributions to the theory of the Riemann zeta-function and the theory of the distribution of primes''~\cite{HL},
Hardy and Littlewood  derived the formula
\be\label{Ram 1}
\begin{split}
\sqrt a F(a)
- \sqrt b F(b)   
 =-\tfrac{1}{2 }\sum_{\r} \frac{\displaystyle\Gamma
(\tfrac12-\tfrac{\r}{2})}{\Z'(\r)} b^{\r-\frac12},
\end{split}
\ee
where $a, b>0$ and $a b=\pi$.   Here the sum  runs over the nontrivial zeros $\r=\b+i\g$ of the zeta function and 
we have assumed they are all simple  (the sum can be modified accordingly if  they are not). 
The formula was suggested to them by some work of Ramanujan. Hardy and Littlewood mentioned
that there  is a way to bracket the terms in the sum over zeros to ensure convergence, but they were not explicit about how to do this.
  Titchmarsh~\cite{T} (see pp. 219-220), however,  proved that the series converges   provided   any two zeros $\r_1, \r_2$ in the sum  for which
\be\label{zero gap 1} 
  |\g_1-\g_2| \leq e^{-A\g_1/\log \g_1} +e^{-A\g_2/\log \g_2}, 
\ee 
with $A$   a sufficiently small positive constant, are grouped together.
In addition, Hardy and Littlewood   proved that for any $\e>0$, the estimate    
\be\label{F bd}
F(b) \ll_\e b^{-\frac12+\epsilon}     
\ee
as $b\to\infty$ is equivalent to the Riemann hypothesis (RH), and they conjectured that, in fact,
$F(b) \ll b^{-\frac12}$.

Several mathematicians have studied various aspects and analogues of $F(b)$ and Ramanujan's formula. For instance, W. Sta\'s~\cite{S1, S2, S3} proved, under various hypotheses,  results of the form
$$
\max_{T^{1-o(1)} \leq b\leq T} |F(b)| \gg T^{-\frac12-o(1)},
$$
for $T$ sufficiently large. A. Dixit~\cite{D1, D2} proved analogues of \eqref{Ram 1} with Dirichlet characters and the insertion of other functions in the sums.
Other results along similar lines may be found in  \cite{Agar, JMS, KRR, RZZ} to cite just a few examples.

Our  purpose here  is to record a few observations about the  finer behavior of $F(b)$ as well as the sum over zeros on the right-hand side of \eqref{Ram 1}  under the assumption of two well-known and widely believed hypotheses. 
We will refer to our first hypothesis as the  {\emph{weak Mertens hypothesis}}  
(WMH).
\begin{hyp 1} \label{con: WMH} 
Let $M(x)=\sum_{n\leq x} \mu(n)$. Then as $X\to\infty$,
\be\label{eq: WMH} 
\int^X_1 \Big( \frac{M(x)}x \Big)^2 dx \ll \log X.
\ee
\end{hyp 1}

We  assume WMH throughout. 
It has the following  consequences:
\begin{itemize}
\item[(A)]  RH,
\item[(B)]   all the zeros $\r$ are simple,
\item[(C)]    $\zeta'(\r)^{-1} =o( |\r|)$,
\item[(D)]  there is a positive constant $A$ such that if $\g<\g'$ are consecutive ordinates of nontrivial zeros of $\zeta(s)$, then
\be\label{zero gap 2}
\g'-\g >  \frac{A}{\g}  \exp\Big(-A \frac{\log \g}{\log\log \g}\Big).
\ee
\end{itemize}
For proofs that WMH implies  (B), (C), and (D), we refer the reader to Titchmarsh~\cite{T}  (Sections 14.29, 14.31).
The proof that WMH implies RH is not in Titchmarsh, but it is short so we provide it here. Set
$$
f(x) =\int^x_1 \frac{M(u)}u  du. 
$$
By the Cauchy-Schwarz inequality and \eqref{eq: WMH},
$$
f(x)^2 \leq x \int^x_1 \Big(\frac{M(u)}u\Big)^2  du \ll   {x \log x}.
$$
Hence $f(x)\ll (x\log x)^\frac12$.
Thus, for $s=\s+it$ with $\s>1$
\be\label{recip zeta}
\begin{split}
\frac{1}{s\zeta(s)} =&\int_1^\infty \frac{M(x)}{x^{s+1}} dx
= \int_1^\infty \frac{d f(x)}{x^{s}} dx
= s \int_1^\infty \frac{f(x)}{x^{s+1}} dx,
\end{split}
\ee
 and it follows that the  last integral in \eqref{recip zeta} is
an analytic function for $\s>1/2$. Thus, $\zeta(s)$ has no zeros in $\s>1/2$. In other words, RH follows.
 

From \eqref{zero gap 2} we see that there are no  zeros with ordinates
$\g_1, \g_2$ large  such that \eqref{zero gap 1} holds. Thus, assuming WMH, 
 \eqref{Ram 1} holds with the sum interpreted as $\lim_{T_\nu\to\infty}\sum_{|\g|\leq T_\nu}$
 for \emph{any} increasing sequence $\{T_\nu\}$. However, on WMH  even more is true -- the series is in fact absolutely convergent. 
 To see this, write
 \be\notag
\sum_{\r} \frac{\displaystyle\Gamma
 (\tfrac12-\tfrac\r{2})}{\Z'(\r)} b^{\r-\frac12}
= \sum_{\g} a(\g)b^{i\g}.
\ee
By Stirling's formula, 
\be\notag
\log \Gamma(s) =(s-\tfrac12)\log s -s +\tfrac{1}{2} \log 2\pi +O(|s|^{-1}),
\ee
where $|s|\to\infty$ in any angle $-\pi+\d<\arg s<\pi-\d$ with $\d>0$. 
Thus
\be\notag
\log |\Gamma  ( \tfrac12-\tfrac{\r }{2} ) | 
= -\tfrac\pi{4} |\g|  -\tfrac{1}{4} \log |\g|+O(1).
\ee
Using this and (C), we find that
\be\label{a bd}
a(\g) =o(|\g|^{3/4} e^{-\pi |\g|/4}).
\ee
Hence, since $N(T) =\sum_{0<\g\leq T} 1 \sim (T/2\pi) \log T$ and the zeros $\r=\frac12+i\g$ are symmetric about the real axis,
we have
\be\notag
\sum_{\g} |a(\g) b^{i\g}|
=O\Big(\sum_{\g}  |\g|^{3/4} e^{-\pi |\g|/4} \Big) \ll 1.
\ee

Returning to \eqref{Ram 1}, we see that 
since  the zeros $\r=\frac12+i\g$ are symmetric about the real axis and $\zeta'(s)$ and $\Gamma(s)$ are real on the real axis, we may  rewrite \eqref{Ram 1} as
\be\notag
\begin{split}
\sqrt a F(a)
- \sqrt b F(b)   
 =- \Re  \sum_{\g>0} a(\g) b^{i\g}.
\end{split}
\ee
Then, since $a b=\pi$ with $a, b>0$, we may replace $a$ by $\pi/b$ and write
\be\label{Ram 3}
\begin{split}
F(b)   
 =\frac{1}{\sqrt b} \Re  \sum_{\g>0} a(\g) b^{i\g}  +\frac{\sqrt\pi}{b}F\Big(\frac\pi{b}\Big),
\end{split}
\ee
 where the sum over $\g$ on the right-hand side  is absolutely convergent under the assumption of WMH.
Since $  \sum_{n=1}^{\infty}  {\mu(n)}n^{-1} =0$, we have 
\be\label{F small arg}\notag
\begin{split}
F\Big( \frac{\pi}{b}\Big) =&\sum_{n=1}^{\infty} \frac{\mu(n)}{n}  (e^{-(\pi/bn)^2}  -1 )  
= \sum_{n=1}^{\infty} \frac{\mu(n)}{n} \sum_{k=1}^{\infty} \frac{(-1)^k (\pi/b n)^{2k}}{k!}\\
=& \sum_{k=1}^{\infty} \frac{(-1)^k (\pi/b)^{2k}}{k!}\sum_{n=1}^{\infty} \frac{\mu(n)}{n^{2k+1}} 
= \sum_{k=1}^{\infty} \frac{(-1)^k (\pi/b)^{2k}}{k! \zeta(2k+1)},
\end{split}
\ee
where the interchange of summations is justified by absolute convergence.
For $b\geq \pi$ it is easily checked that the absolute value of the terms of this alternating series are decreasing, so for any integer $K\geq 1$, we see that
\be\notag
F\Big( \frac{\pi}{b}\Big) = \sum_{k=1}^{K} \frac{(-1)^k (\pi/b)^{2k}}{k! \zeta(2k+1)} +E_{K+1}(b),
\ee
where
\be\notag
|E_{K+1}(b)|\leq \frac{(\pi/b)^{2K+2}}{(K+1)!} .
\ee
Inserting this into \eqref{Ram 3}, we now find that if WMH is true and $b\geq \pi$, then
\be\label{Ram 4}
\begin{split}
F(b)   
 =\frac{1}{\sqrt b} \Re  \sum_{\g>0} a(\g) b^{i\g}  
 +\frac{1}{\sqrt\pi} \sum_{k=1}^{K} \frac{(-1)^k (\pi/b)^{2k+1}}{k! \zeta(2k+1)}
    +\frac{\sqrt \pi}{b}E_{K+1}(b).
\end{split}
\ee
We will use this for the calculations in Section~\ref{Calculations}. 
However, even the cruder estimate 
 \be\label{Ram 5}
\begin{split}
  F(b)   
 =  \Re \;\frac{1}{ \sqrt b }\sum_{\g>0}\; a(\g) b^{i\g}
+O(b^{- 3})
\end{split}
\ee
immediately  leads to the following theorem.
\begin{thm}\label{thm 1} Assume WMH. Then for $b\geq \pi$ we have
\be\notag
|F(b)| \leq \frac{C}{\sqrt b} +O(b^{-3}),
\ee 
 where
\be\notag 
C= \sum_{\g>0} |a(\g)| = \sum_{\g>0} \bigg|\frac{ \Gamma
 (\frac{1}{4}-i\frac{\g}{2} )}{\Z'(\frac12+i\g)}\bigg|.
\ee
\end{thm}
  
To  analyze the sum over $\g$ in \eqref{Ram 4} and \eqref{Ram 5}, 
we assume, in  addition to WMH,  the following \emph{linear independence hypothesis} (LI).
\begin{hyp 2}
The positive ordinates $\g$ of the zeros of the zeta function are linearly independent over the rationals. 
\end{hyp 2}

To use this  we first assume the  $\g>0$ have been ordered   as $\g_1, \g_2, \g_3 \ldots,$  in such a way  that  $|a(\g_1)|\geq |a(\g_2)|\geq |a(\g_3)|\geq\cdots$.
Then
$$
\sum_{\g>0}\; \frac{ \Gamma(\tfrac14-i\tfrac\g2 )}{\Z'(\frac12+i\g)} b^{i\g}
=\sum_{\g>0}\; a(\g)  b^{i \g}
=\sum_{n=1}^{\infty}\;   a(\g_n)  e^{i \g_n \log b}.
$$
 LI  implies that as $b$ varies over $[\pi, \infty)$, this sum   is dense in the set of complex numbers
$$
\mathscr A =\Big\{ \sum_{n=1}^{\infty}\;   |a(\g_n)|  e^{i \theta_n }\,:\, \theta_n\in[0,1), n=1, 2, 3, \ldots \Big\}.
$$ 
This set, being a ``sum'' of  circles centered at the origin, is, as is well known, either  a closed annulus
or a closed disk according to the following criteria:

 
\noindent 1) If 
 $ \displaystyle  |a(\g_1)|>\sum_{n= 2}^\infty |a(\g_n)| $,
then $\mathscr A$ is a closed annulus centered at the origin with  outer 
radius 
$$
C=\sum_{n= 1}^\infty |a(\g_n)|
$$ 
and  inner radius 
$$c= |a(\g_1)|-\sum_{n= 2}^\infty |a(\g_n)|.$$
\noindent 2) If 
$\displaystyle 
|a(\g_1)|\leq \sum_{n= 2}^\infty |a(\g_n)|,$ then $\mathscr A$ is a closed disk centered at the origin of radius 
$$
C=\sum_{n= 1}^\infty |a(\g_n)|.
$$
 
%

In either of these two cases,  the real parts of the complex numbers 
$\sum_{n=1}^{\infty}\;   |a(\g_n)| e^{i \theta_n}$ in $\mathscr A$ fill out the  
  interval $[-C, C]$. As the sum $\sum_{n=1}^{\infty}\;   a(\g_n)  e^{i \g_n \log b}$ is dense in $\mathscr A$ 
(assuming LI),  this and \eqref{Ram 5} give the following result.
\begin{thm} Assume WMH and LI.
Then $\sqrt b F(b)$ is dense in $[-C, C]$ and, in particular, we have
$$
\liminf_{b\to\infty}\sqrt b F(b) =-C
\qquad\qquad \hbox{and}  \qquad\qquad \limsup_{b\to\infty} \sqrt b F(b) =C .
$$
\end{thm}
%

For $N$ a large positive integer, let 
$$
\mathscr A_N =\Big\{ \sum_{n=1}^{N}\;   |a(\g_n)|  e^{i \theta_n }\,:\, \theta_n\in[0,1), n=1, 2, 3, \ldots, N\Big\},
$$ 
which again is either an annulus or disk centered at the origin.
By the reasoning above, if one assumes  LI, the curve $f_N(b)=\sum_{n=1}^{N}\;   a(\g_n)  e^{i \g_n \log b}$ is   dense in $\mathscr A_N$. By the Kronecker-Weyl theorem, it is also uniformly distributed in  $\mathscr A_N$. Thus, the distribution function of the curve $\Re f_N(b)$ as $b\to\infty$ tends to the distribution function of the $x$ coordinate of points $(x, y)$ in the annulus or disk $\mathscr A_N$. Since $\sum_{n=1}^{\infty}\;   a(\g_n)  e^{i \g_n \log b}$ is absolutely convergent, the same is true for the real part of this  series but with $\mathscr A$ in place of 
$\mathscr A_N$.
Moreover, by \eqref{Ram 5},  
$$\Big|\sqrt{b}F(b)- \Re \; \sum_{n=1}^{\infty}\;   a(\g_n)  e^{i \g_n \log b}\Big| \ll b^{-\frac52}.$$
Thus, as $b\to \infty$, the probability distribution function of  $\sqrt{b}F(b)$ tends to the distribution function of the $x$ coordinate of points $(x, y)$
 in either the annulus centered at the origin with inner radius $c$ and outer radius $C$,  or  the disk centered at the origin of radius $C$. 
Depending on whether the set $\mathscr A$ is an annulus or a disk, we therefore have  the following probability density function for $\sqrt b F(b)$.

\begin{thm} Assume WMH and LI. Let $c$ and $C$ be as above, let
$$
\mathscr A =\Big\{ \sum_{n=1}^{\infty}\;   |a(\g_n)|  e^{i \theta_n }\,:\, \theta_n\in[0,1), n=1, 2, 3, \ldots \Big\},
$$ 
and let $p(x)$ be the probability density function   of $\sqrt b F(b)$ for $b$ large.
If $\mathscr A$ is an annulus with inner radius $c$ and outer radius $C$, then
\be\notag
p(x)=
\begin{cases}
 \ds 0 &\quad\text{if} \quad  C\leq |x|, \\
\ds \frac{2 \sqrt{C^2-x^2}}{\pi (C^2-c^2)}  &\quad\text{if} \quad c\leq |x|\leq C,\\
\ds \frac{ 2 (\sqrt{C^2-x^2} - \sqrt{c^2-x^2} )}{ \pi (C^2-c^2) } &\quad\text{if} \quad   |x|\leq c.
\end{cases}
\ee
If $\mathscr A$ is a disk of radius $C$, then
\be\notag
p(x)=
\begin{cases}
\ds 0 &\quad\text{if} \quad  C\leq |x|, \\
\ds \frac{2 \sqrt{C^2-x^2}}{\pi C^2}   &\quad\text{if} \quad  |x|\leq C.
 \end{cases}
\ee
\end{thm}


It seems difficult to prove, even under the strong assumptions of WMH and LI, whether $\mathscr A$ is an annulus or disk, but we believe it to be an annulus. At issue is determining the relative size of the two quantities
$$
 |a(\g_1)|  \qquad \hbox{and} \qquad \sum_{n= 2}^\infty |a(\g_n)|,
 $$
 where
 $$
 a(\g) =  \frac{ \Gamma(\tfrac14-i\tfrac\g2 )}{\Z'(\tfrac12+i\g)}.
 $$
There are two sources of difficulty in settling this question. One is that,
although the size of $\Gamma$ is well understood,  the   bound 
$\zeta'(\r)^{-1}=o(|\r|)$ from C) is not explicit enough; what would suffice  is an estimate of the type
$|\zeta'(\r)^{-1}|\leq B|\r|$ for all $\g>0$ with $B$ an explicit constant, or even $|\zeta'(\r)^{-1}|\leq B |\g|^d $ with $d>1$, and $d$ and $B$ both explicit. The other difficulty, which is  related to the first, is  that we do not know which $\g$ should be $\g_1$, that is, which $\g$ maximizes $|a(\g)|$. (Note that if $|a(\g)|$ is maximal for more than one $\g$, then $\mathscr A$ is a disk.) 
However, if a constant $B$ as above exists that is not enormous, the fast exponential decay from the gamma function in $a(\g)$ suggests that the drop off between terms for  successive $\g$'s is large, and this suggests that $a(\g_1)$ (with $\g_1=\g$) is much larger than $\sum_{n= 2}^\infty |a(\g_n)|$.  In Section~\ref{Calculations} we present the outcome of a limited number of calculations that suggest possible approximate values of $c$ and $C$ and  we present several graphs of $\sqrt b F(b)$.

We next prove a formula for the second moment of $F$.
\begin{thm}\label{mean square} Assume WMH. Then 
\be\label{eq: mean square}
\int_1^X F(x)^2 dx = A \log X +O(1)
\ee
as $X\to\infty$, where
$$A=\tfrac12 \sum_{\g>0}\; |a(\g)|^2.$$
\end{thm}

\noindent{\emph{Remark.}} Note that $A>0$.
\begin{proof}

Writing 
$$
S = \sum_{\g>0}\; a(\g) x^{i\g} ,
$$
we find by \eqref{Ram 3} that
 \be\notag
\begin{split}
 \int_1^X F(x)^2 dx  
 = & \int_1^X \Big(\frac{\Re S}{ \sqrt x }+O(x^{- 3}) \Big)^2 dx\\
 = & \int_1^X \big(\Re  S +O(x^{- 5/2}) \big)^2 \frac{dx}{x}\\
 = & \int_1^X \big( (\Re  S )^2 +O(|S|x^{- 5/2}) +O(x^{-5}) \big) \frac{dx}{x}.
\end{split}
\ee
Since the series defining $S$ is absolutely convergent, the last two terms of the integrand contribute  
$O(1)$. Thus,
\be\label{2nd mmt}
 \begin{split}
 \int_1^X F(x)^2 dx  = &\int_1^X (\Re S )^2 \frac{dx}{x} +O(1) \\
 =&\frac14  \int_1^X ( S^2 +2|S|^2 +\overline{S}^2 ) \frac{dx}{x} +O(1) . 
\end{split}
 \ee
Again, by absolute convergence of the sum defining $S$, we have 
 \be\notag
 \begin{split}
 \int_1^X S^2 \frac{dx}{x} = &\int_1^X \sum_{\g, \g' >0}\; a(\g) a(\g') x^{i(\g+\g')}  \frac{dx}{x} 
 = \sum_{\g, \g' >0}\; a(\g) a(\g') \int_1^X  x^{i(\g+\g')-1}   dx \\
  = &\sum_{\g, \g' >0}\; a(\g) a(\g')    \frac{X^{i(\g+\g')}-1 }{i(\g+\g')}   
   \ll \Big( \sum_{\g>0}\; |a(\g)| \Big)^2 \ll 1.
 \end{split}
 \ee
 Similarly, $\ds \int_1^X {\overline{S}}^2 \frac{dx}{x} \ll 1$. Finally,
 \be\notag
 \begin{split}
 \int_1^X  |S|^2   \frac{dx}{x}  
  =& \sum_{\g, \g' >0}\; a(\g) \overline{a(\g') } \int_1^X  x^{i(\g-\g')-1}   dx \\
  = &\log X \sum_{\g>0}\; |a(\g)|^2 +  \sum_{\substack{\g, \g' >0\\ \g\neq \g'}}\; a(\g)  \overline{a(\g') }   \frac{X^{i(\g+\g')}-1 }{i(\g-\g')} \\
  = & \log X  \sum_{\g>0}\; |a(\g)|^2+ O\bigg( \sum_{\substack{\g, \g' >0\\ \g\neq \g'}}\; |a(\g)  \overline{a(\g') }  | \min\Big( \log X,   \frac{1}{|\g-\g'|}  \Big) \bigg).
 \end{split}
 \ee
By (D), for any  $\e>0$ we have $|\g-\g'|^{-1} \ll \g^{1+\e}$. Thus,  by \eqref{a bd}, we find that the $O$-term is
 \be\notag
 \begin{split}
\ll  & \sum_{ \g' >0}  |a(\g')|\sum_{ \g<\g'}  | a(\g) | \g^{1+\e}
\ll   \sum_{ \g' >0}  |a(\g')| \sum_{ \g<\g'}   \g^{7/4+\e} e^{-\pi |\g|/4}\\
\ll & \sum_{ \g' >0}  |a(\g')| \ll 1.
 \end{split}
 \ee
 Hence
 \be\notag
 \begin{split}
 \int_1^X  |S|^2   \frac{dx}{x}  
  = & \log X  \sum_{\g>0}\; |a(\g)|^2+ O(1).
 \end{split}
 \ee
 Combining our estimates together in \eqref{2nd mmt}, we obtain
 \be\notag
 \int_1^X F(x)^2 dx =\log X  \Big(\tfrac12 \sum_{\g>0}\; |a(\g)|^2\Big) +O(1) . 
 \ee
 \end{proof}
\noindent\emph{Remark.} \
One can show that if a weak version of \eqref{eq: mean square} holds,  namely,
$$
\int_1^X F(x)^2 dx \ll \log X,
$$
then   (A) and (B) as well as  the following analogue of (C) follow:
\begin{itemize}
\item[(C*)]   $\zeta'(\r)^{-1} \ll e^{c|\g|}$ for some positive constant $c$.
\end{itemize}
These can be proved along the lines of the proofs that  (A)-(C) follow from WMH.\\

\section{Riesz's function}

Analogues of the results above may easily be extended to M. Riesz's function~\cite{R}
$$
P(x)=\sum_{n=1}^\infty\dfrac{\mu(n)}{n^2}e^{-x/n^2},
$$
which is similar to $F(x)$ and was introduced   around the same time as Hardy and Littlewood's 
work on Ramanujan's formula.
Note that $P(x)$ has $n^2$ rather than $n$ in the denominator and $x$ rather than $x^2$ in the exponential.
 Agarwal, Garg, and Maji~\cite{Agar} recently generalized this to a one parameter family of functions
$$
P_k(x)=\sum_{n=1}^\infty\dfrac{\mu(n)}{n^k}e^{-x/n^2},
$$
where $k\geq 1$ is a fixed real number. Note that $F(x)=P_1(x^2)$ and $P(x)=P_2(x)$.
They then proved the following analogue of \eqref{Ram 1} (see their Theorem 1.1): 
\begin{align} \label{Agarwal form 1}
P_k(x)= \Gamma(\tfrac{k}{2}) x^{-\frac{k}{2}}\sum_{n=1}^{\infty}\dfrac{\mu(n)}{n}\ {_1F_1}\big(\tfrac k2;\tfrac{1}{2};-\tfrac{\pi^2}{n^2x}\big) + \tfrac{1}{2}\sum _{\rho}a_k(\g)x^{-\frac{k-\r}{2}}.
\end{align}
Here ${_1F_1}(\frac{k}{2};\hh; z)$ is the generalized hypergeometric series, 
$$
a_k(\r)=\frac{\Gamma(\frac{k-\rho}{2})}{\zeta'(\rho)},
$$
the zeros $\r$ are all assumed to be simple, and  any two zeros $\r_1$ and $\r_2$ in  the series on the right in \eqref{Agarwal form 1} are grouped together if they satisfy the inequality \eqref{zero gap 1}.
 They  used this to show that for any fixed real number $k\geq 1$ and any  $\epsilon>0$, the Riemann hypothesis is equivalent to 
\be\notag
P_k(x) \ll_\epsilon  x^{-\frac k2+\frac14+\epsilon}
\ee
as $x\to\infty$ (similarly to \eqref{F bd}).

  Assuming WMH and  using \eqref{Agarwal form 1}, we may easily prove a version of \eqref{Ram 5} for $P_k(x)$.
 First note, as before, that from WMH
  it follows that RH holds, all the zeros $\rho$ of $\zeta(s)$ are simple, and $|\zeta(\rho)^{-1}|=o(|\rho|)$. Also, by Stirling's formula, we have 
$$
\log\left|\Gamma\left(\tfrac{k}{2}-\tfrac{\rho}{2}\right)\right|
= - \tfrac{\pi}{4} |\gamma|+ (\tfrac{k}{2}-\tfrac{3}{4}) {\log|\gamma|} + O(1).
$$
Thus,
\begin{align*}
\sum_{\r} |a_k(\r)| 
	\ll &   \sum _{\gamma}|\gamma|^{{\frac{k}{2}+\frac{1}{4}}}e^{-\pi|\gamma|/4}\ll 1.
\end{align*}
Hence, the series  
$$
 \tfrac{1}{2}\sum _{\rho}a_k(\rho)x^{-\frac{k-\rho}{2}} 
 = \tfrac{1}{2}x^{-\frac{k}{2}+\frac{1}{4}} \sum _{\rho}a_k(\rho)x^{ i\g/{2}} 
 =x^{-\frac{k}{2}+\frac{1}{4}} \Re \ \sum _{\g>0}a_k(\rho)x^{ i\g/{2}} 
$$
on the right-hand side of \eqref{Agarwal form 1} converges absolutely.
Next,  for    $z$ complex and bounded, we have 
 \be\notag
 \begin{split}
 {_1F_1}(\tfrac{k}{2},\hh,z) =\sum_{j=0}^{\infty} \frac{  \Gamma(\tfrac k2+j) \Gamma(\tfrac 12)  }{ \Gamma(\tfrac k2) \Gamma(\tfrac12+j) j!} z^j
 =1+O(|z|).
 \end{split}
 \ee
Thus, the first term on the right-hand side of \eqref{Agarwal form 1} equals
$$
\Gamma(\tfrac{k}{2}) x^{-\frac{k}{2}}\sum_{n=1}^{\infty}\frac{\mu(n)}{n} (1+O(n^{-2}x^{-1} ) ) 
\ll x^{-\frac k2-1},
$$
since $\sum_{n}\mu(n)n^{-1}=0$.

Using these estimates and observations with \eqref{Agarwal form 1}, we arrive at
$$
P_k(x)= x^{-\frac{k}{2}+\frac{1}{4}} \ \Re \ \sum _{\g>0}a_k(\rho)x^{ i\g/{2}}+ O({x^{-\frac{k}{2}-1}}).
$$
With this formula as a starting point, we may easily prove analogues of Theorems 1-4  for $P_k(x)$. In the case of Theorem 4,
we obtain an asymptotic  formula for 
$$
\int_1^X P_k(x)^2 x^{k-\frac32} dx.
$$ 
\smallskip

\section{Calculations}\label{Calculations}
We mentioned in Section~\ref{Ram form} that we believe $\mathscr A$ to be an annulus.
 In this final section we briefly report  the results of  calculations of  a number of $|a(\g)|$'s, 
 and use these to approximate the values of the inner and out radii, $c$ and $C$, of the annulus $\mathscr A$.   We also provide several graphs of $\sqrt b F(b)$. We have used Mathematica 
for these calculations and to generate our graphs. 

First, here is a table of values of $|a(\g)|$ for the first ten ordinates $\g>0$.\\
 
\hskip2in
\begin{tabular}{|l | c | }
\hline
$\gamma_n$& $|a(\gamma_n)|$ \\ \hline 
$\gamma_1$ &$2.9255\ldots\cdot10^{-5}$ \\ \hline 
$\gamma_2$ & $8.2702\ldots\cdot10^{-8}$ \\ \hline
$\gamma_3$ & $2.8609\ldots\cdot10^{-9}$ \\ \hline
$\gamma_4$ & $4.0789\ldots\cdot10^{-11}$ \\ \hline
$\gamma_5$ & $5.2534\ldots\cdot10^{-12}$ \\ \hline
$\gamma_6$ & $9.4006\ldots\cdot10^{-14}$ \\ \hline
$\gamma_7$ & $8.7272\ldots\cdot10^{-15}$ \\ \hline
$\gamma_8$ & $1.0550\ldots\cdot10^{-15}$ \\ \hline
$\gamma_9$ & $3.0507\ldots\cdot10^{-17}$ \\ \hline
$\gamma_{10}$ & $8.3287\ldots\cdot10^{-18}$ \\ \hline
\end{tabular}

Notice that, for the most part, these terms are quickly decreasing.
If we sum them to approximate $C$, the outer radius of $\mathscr A$, we obtain the value  
$C\approx 0.0000293414$. To approximate $c$ we subtract the sum of the last nine values from $|a(\g_1)|$
and obtain $c\approx 0.0000291702$. Interestingly, performing the same calculations with the first $500$ ordinates 
$\g$ gives exactly the same values for $C$ and $c$ up to ten significant figures. This suggests (but, of course, does not prove) that $\mathscr A$ really is an annulus rather than a disk.

We conclude with several graphs of $\sqrt b F(b)$ for various ranges of $b$ from the formula \eqref{Ram 4} using the first $50$ ordinates $\g$ and the sum over $k$ with 
$ K=50$ and ignoring  the error term  $E_{51}(b)$. 
Although our estimate for $E_{K+1}$ in \eqref{Ram 4} was for $b\geq \pi,$ it is not difficult to check that $E_{51}(b)$ is quite small even when $1\leq b\leq \pi$. Thus, Figure 1 is accurate for this range of $b$ as well.

For some related graphs see Paris~\cite{P}.

\begin{figure}[ht]\label{function 1}
	\centering 
	\includegraphics[width=3.8in]{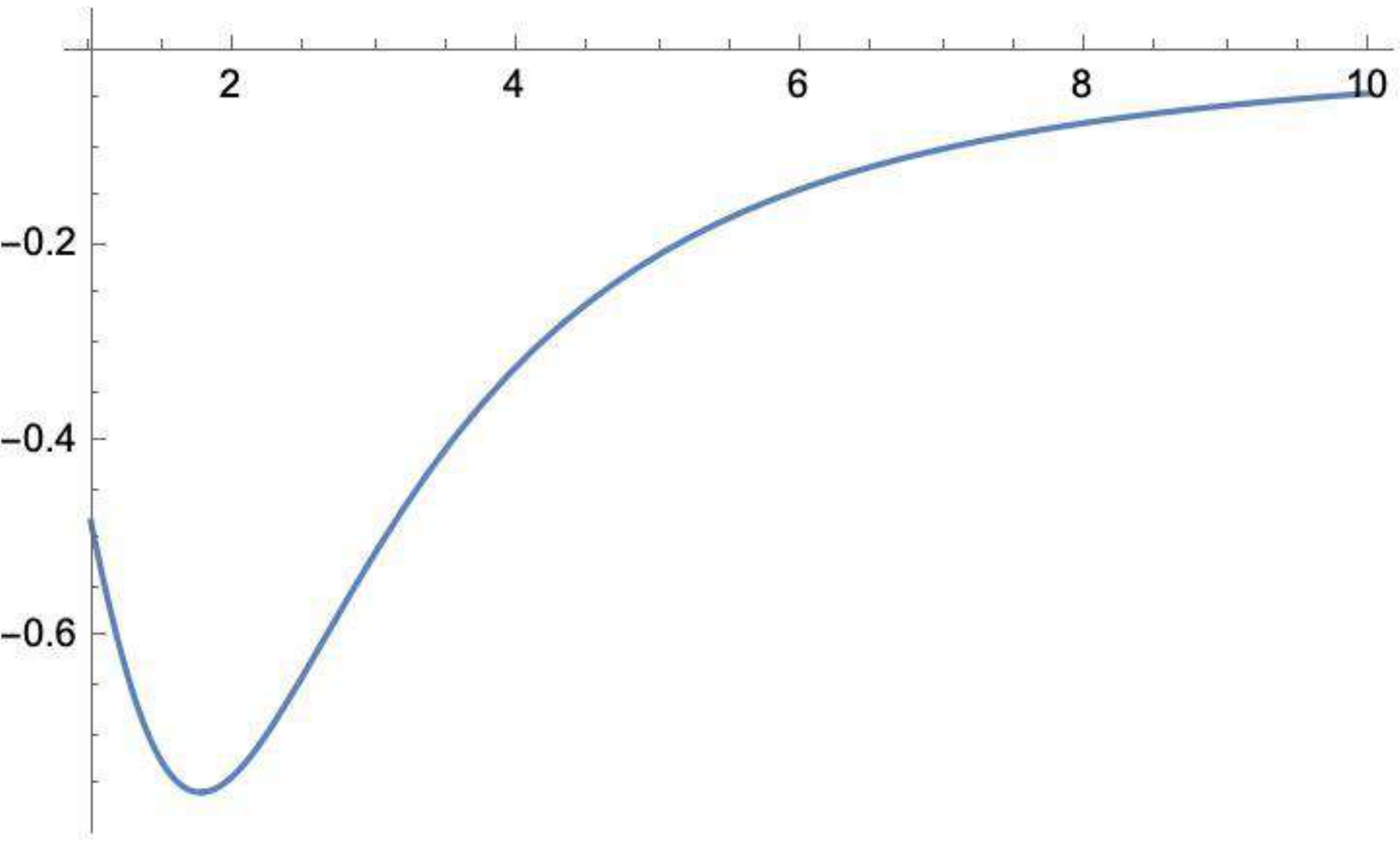}
	\vspace{-0.4cm}\caption{Graph of  $\sqrt b F(b)$ for $1\le b\le 10$.} \label{function1}
\end{figure}

\medskip

\begin{figure}[ht]\label{function 2}
	\centering 
	\includegraphics[width=4.2in]{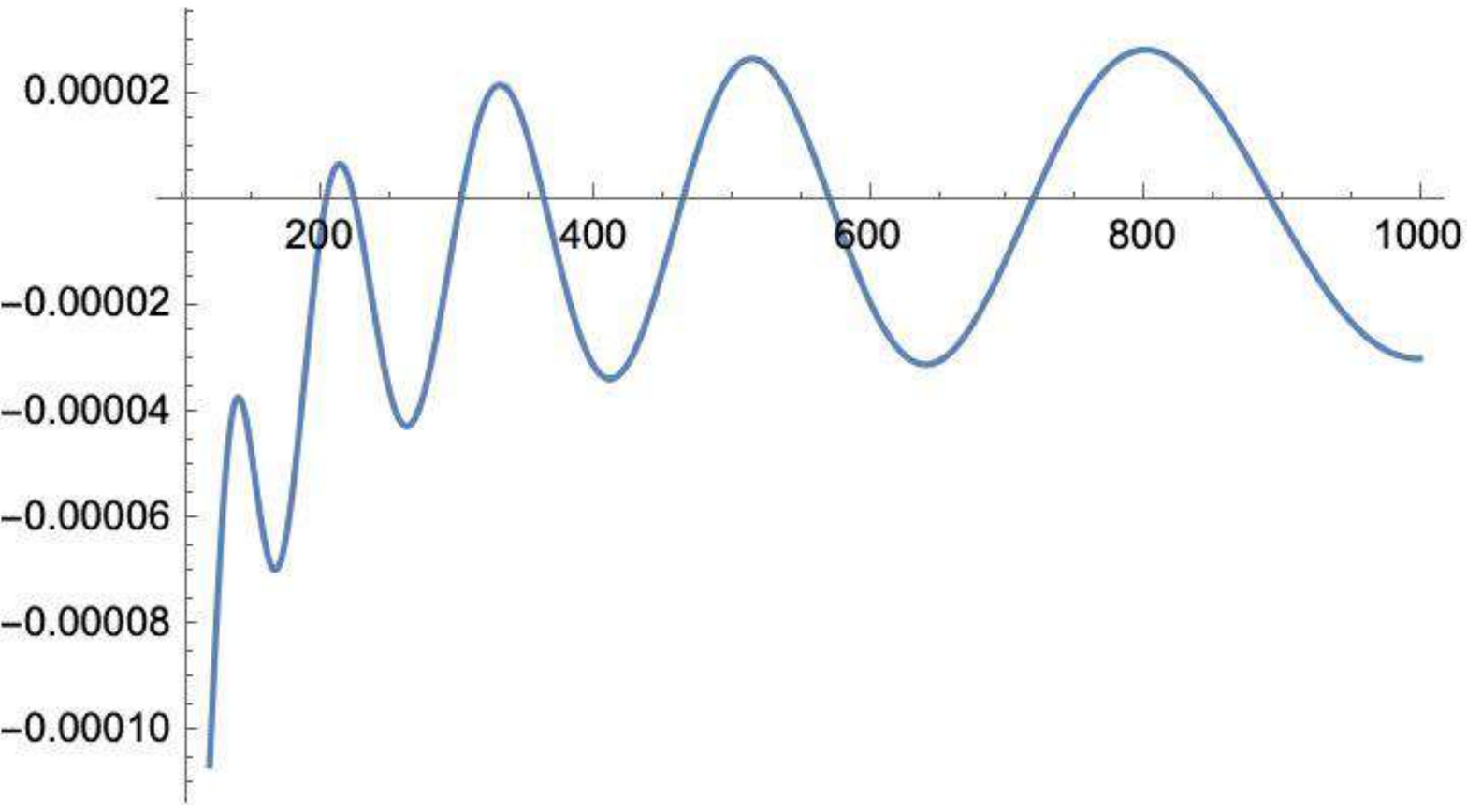} \vspace{-0.4cm}
	\caption{Graph of  $\sqrt b F(b)$ for $100\le b\le 1000$.}\label{function 2}
\end{figure}

\begin{figure}[ht]\label{function 3}
	\centering 
	\includegraphics[width=4.2in]{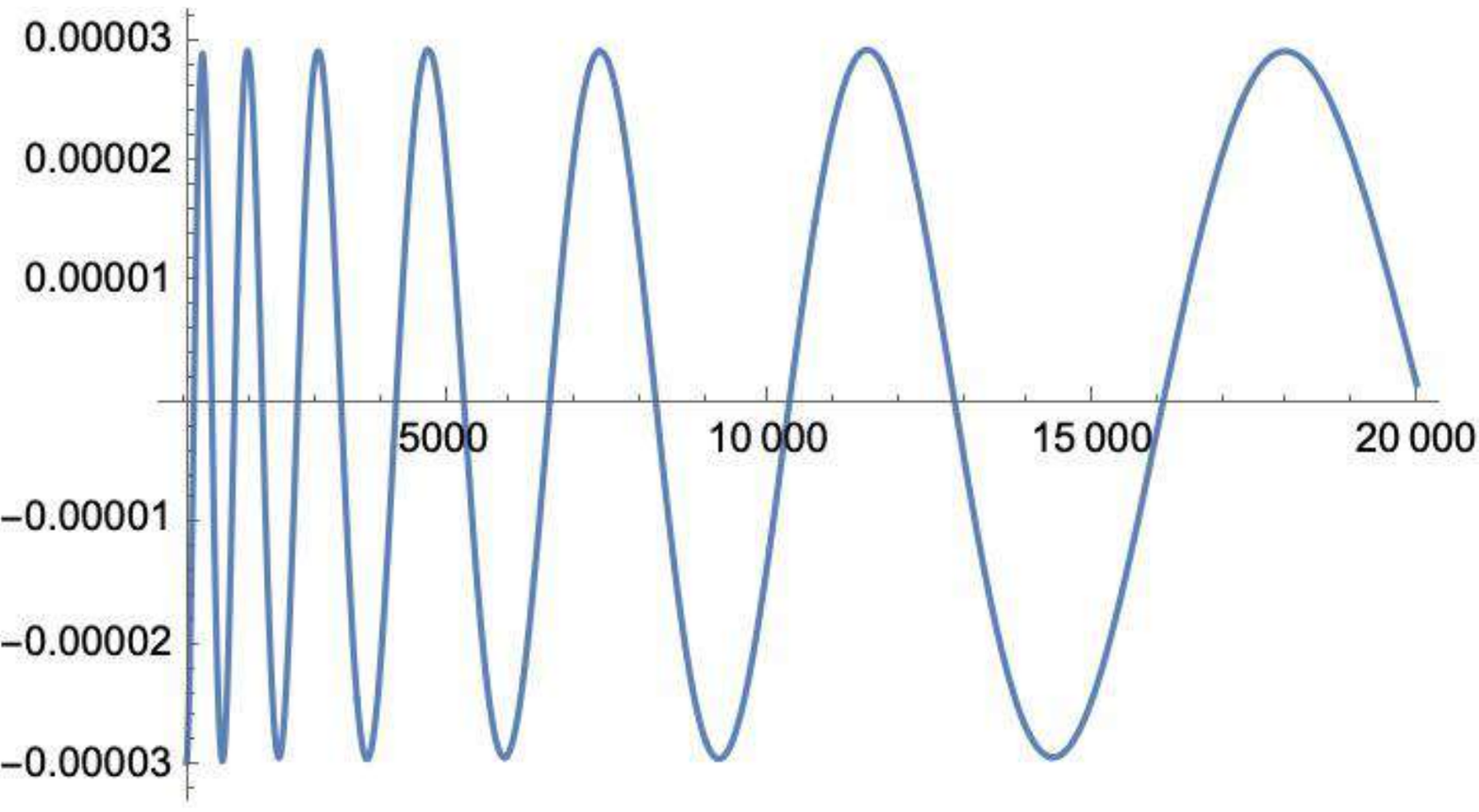}\vspace{-0.4cm}
	\caption{Graph of  $\sqrt b F(b)$ for $1000\le b\le 20000$.}\label{function 3}
\end{figure}

\newpage

\bigskip

\section{acknowledgements}
We thank Dr. Shashank Chorge for carrying out a preliminary version of the calculations leading to our approximate values  of 
$C$ and $c$. This project was begun when Andr\'es Chirre was a Visiting Assistant Professor at the University of Rochester. He thanks the University of Rochester for its hospitality and support.


\end{document}